\documentclass[preprint,12pt]{article}
\usepackage{tikz}
\begin{document}

\newtheorem {theorem} {Theorem}
\newtheorem {lemma} {Lemma }

\title{There is no Perfect  Cuboid }

\author{Ivor Lloyd  \\  \textit{(Greenfields, 32A Brook Hill, Woodstock , OX20 1JE ,UK \footnote{Email. ivor.lloyd@psy.ox.ac.uk})}}

\maketitle

\begin{abstract}

A perfect cuboid is formed when an Euler brick whose edges and face diagonals are all integers also has an integer internal diagonal. It is known that if a perfect cuboid exists the internal diagonal is odd. No perfect  cuboid has been found. This simple proof shows that the internal diagonal of an Euler brick cannot be an odd integer.
\end{abstract}

\section {Introduction}

A perfect cuboid is formed when an Euler brick whose edges and face diagonals are all integers also has an integer internal diagonal.  The problem is to find a perfect cuboid or prove that one cannot exist.   Sharipov[2]  states that this problem was first mentioned by Paul Halcke[1] some 300 years ago. Although many Euler bricks have been shown to exist, no perfect  cuboid has yet been found. In this simple proof we show that the internal diagonal of the Euler brick cannot be an odd integer and thus there is no perfect cuboid.

It is known that two of the edges of the Euler brick are even and the other odd[3]. Let the the edges of the cuboid be  integers \begin {math} 2^{n}a,b,2^{n}c \end {math}, the face diagonals be \begin {math} 2^{n}d,e,f \end {math}  and the internal diagonal be \begin {math} g\end {math} .  \begin {math} n \end {math} is a positive integer and \begin {math} n>=1 \end {math}. Then the three equations relating the edges and diagonals respectively are:

\begin {equation}(2^nd)^2=(2^na)^2+(2^nb)^2 \end {equation}
\begin {equation} e^2=c^2+(2^nb)^2\end {equation}
\begin {equation} f^2=c^2+(2^na)^2 \end {equation}
Where all variables are odd integers.
The cuboid is perfect if:
\begin {equation} g^2=(2^na)^2+(2^nb)^2+c^2 \end {equation}
The equation for \begin {math} g \end {math} can also be written as:
\begin {equation}2g^2= (2^nd)^2+e^2+f^2 \end {equation}
\section {The proof}

\begin {theorem}  If all the edges and face diagonals in the Euler brick are positive integers then the internal diagonal is not integer.\end {theorem}
\begin {lemma} If  a  perfect  cuboid exists, the internal diagonal of the Euler brick is odd. \end {lemma}
This lemma has been proved by Roberts[3]

In equation (5) let \begin {math} e=2x+1 \end {math} and \begin {math} f=2y+1 \end {math}, where \begin {math} x \end {math} and \begin {math} y \end {math} are positive integers.
\begin {equation} 2g^2=(2^{n}d)^2+(2x+1)^2+(2y+1)^2 \end {equation}
Expanding the brackets and simplifying
\begin {equation} 2g^2=2^{2n}d^2+4x^2+4x+1+4y^2+4y+1 \end {equation}
\begin {equation} g^2=2^{2n-1}d^2+2x^2+2x+2y^2+2y+1 \end {equation}
\begin {equation} g^2=2[2^{2n-2}d^2+x^2+x+y^2+y]+1 \end {equation}
\begin {equation} g^2=2[2^{2n-2}d^2+x(x+1)+y(y+1)]+1 \end {equation}

In equation (10) \begin {math}x(x+1) \end {math} and  \begin {math}y(y+1) \end {math} are even.  If \begin {math} n =1\end {math}, then \begin {math} [2^{2n-2}d^2+x(x+1)+y(y+1)] \end {math} is odd, so the RHS is not a square and \begin{math} g \end {math} is not integer. If \begin {math} n>1 \end {math} the RHS is even. Since \begin {math} g \end {math} is odd, this is a contradiction.

\begin {math} g \end {math} cannot be an odd integer and Theroem 1 is proved.

\end {document}